\documentclass[final,a4paper]{amsart}

\usepackage{amsmath,amssymb,amsthm,amsfonts}
\usepackage{epsfig}
\usepackage{enumitem}

\usepackage[color,notcite]{showkeys}
\definecolor{labelkey}{rgb}{0.6,0,1}

\definecolor{refkey}{rgb}{0.6,0,0} 

\newcommand{\R}{\mathbb{R}}
\renewcommand{\S}{\mathbb{S}^N}

\newcommand{\eps}{\varepsilon}
\newcommand{\un}{\mathbf{1}}

\def\ab{{\lambda, \gamma}}
\def\xb{\bar{x}}
\def\yb{\bar{y}}

\def\tr{{\rm Tr}}

\newtheorem{defi}{Definition}

\newtheorem{thm}{Theorem}
\newtheorem{cor}{Corollary}

\theoremstyle{remark}
\newtheorem{rem}{Remark}
\newtheorem*{rem*}{Remark}

\theoremstyle{remark}

\theoremstyle{remark}

\begin{document}

\title[H\"older continuity and non-linear elliptic non-local equations]%
{H\"older continuity of solutions of 
second-order non-linear elliptic integro-differential equations}

\author{Guy Barles}
\address{Laboratoire de Math\'ematiques et Physique Th\'eorique,
CNRS UMR 6083, F\'ed\'eration Denis Poisson,
Universit\'e Fran\c{c}ois Rabelais, Parc de Grandmont,
37200 Tours, France} 
\email{barles@lmpt.univ-tours.fr}

\author{Emmanuel Chasseigne}
\address{Laboratoire de Math\'ematiques et Physique Th\'eorique,
CNRS UMR 6083, F\'ed\'eration Denis Poisson,
Universit\'e Fran\c{c}ois Rabelais, Parc de Grandmont,
37200 Tours, France}
\email{emmanuel.chasseigne@lmpt.univ-tours.fr}
 
\author{Cyril Imbert}
\address{CEREMADE, UMR CNRS 7534, Universit\'e Paris-Dauphine, Place de
  Lattre de Tassigny, 75775 Paris cedex 16, France}
\email{imbert@ceremade.dauphine.fr}

\maketitle

\begin{quote} \footnotesize
\noindent \textsc{Abstract.} 
This paper is concerned with H\"older regularity of viscosity solutions of second-order, fully non-linear elliptic 
integro-differential equations. Our results rely on two key ingredients: first we assume that, at each point 
of the domain, either the equation is strictly elliptic in the classical fully non-linear sense, or (and this 
is the most original part of our work)  the equation is strictly elliptic in a non-local non-linear sense we make 
precise. Next we impose some regularity and growth conditions on the equation.  These results are concerned with 
a large class of  integro-differential operators whose singular measures depend on $x$ and also a large class of
 equations, including  Bellman-Isaacs Equations.
\end{quote}
\vspace{5mm}

\noindent
\textbf{Keywords:} H\"older regularity, 
integro-differential equations, L\'evy operators, 
general non-local operators, viscosity solutions
\medskip

\noindent  \textbf{Mathematics Subject Classification:} 35D99, 35J60, 35B05, 47G20

\section*{Introduction}

The aim of this paper is to show that viscosity solutions of fully non-linear elliptic integro-differential 
equations are H\"older continuous under general suitable strict ellipticity 
and regularity/growth conditions on the equations. We also obtain explicit $C^{0, \alpha}$ estimates 
{in terms of the space dimension, the non-linearity, the singular measure and the (local) $L^\infty$-bound 
of the solution}.

To be more specific, we describe the general framework we consider. We are interested in equations of the type
\begin{equation}\label{eq:pide}
F(x,u,Du,D^2 u, \mathcal{I} [x,u])=0 \quad \mbox{ in } \Omega
\end{equation}
where $\Omega$ is a domain of $\R ^N$ (not necessarily bounded) and $\mathcal{I} [x,u]$ is an integro-differential 
operator. The function $u$ is real-valued and $Du, D^2 u$ stand respectively for its gradient and Hessian matrix. The non-linearity 
$F$ is a (continuous) function which is degenerate elliptic: in this context, this means that $F$ is 
non-increasing with respect to                 its last two variables (see below for a precise statement and 
\cite{bi07} for further details). 
The integro-differential operators we will consider in the present paper are
defined as follows
\begin{equation}\label{op:id}
\mathcal{I} [x,u] = \int_{\R^N} ( u (x+z ) - u (x) 
- D u (x) \cdot z \un_B (z) ) \mu_x (dz)
\end{equation}
where  $\un_B$ denotes the indicator function of the unit ball $B$ and
$\{ \mu_x \}_{x \in \R ^N}$ is a family of L\'evy measures, 
 \textit{i.e.}  { non-negative, possibly singular, Borel measures on $\R ^N$} such that
\begin{equation}\label{cond:levymeasure}
\int_{\R^N} \min ( |z|^2, 1)\;  \mu_x (dz) < + \infty.
\end{equation}
We point out that the solution $u$ has to be given in the whole space
$\R^N$ even if \eqref{eq:pide} is satisfied on a domain $\Omega$; see
Section~\ref{sec:visc} for further details. 

Important ``special cases'' of operators of the 
form~\eqref{op:id} are L\'evy-It\^o operators, namely
\begin{equation}\label{op:li}
\mathcal{I}_{LI} [x,u]= \int_{\R^N} ( u (x+j (x,z) ) 
- u (x) - D u (x) \cdot j(x,z) \un_B (z) ) \mu (dz)
\end{equation}
where $\mu$ is a L\'evy measure
and $j (x,z)$ is the size of the jumps at $x$. 
For the operator to be well-defined, one usually assumes 
$$
|j(x,z) |\le C_0|z| \quad\hbox{for some constant $C_0$ and for any $x\in \Omega$, $z \in \R^N$}. 
$$
Model examples for a L\'evy measure $\mu$ and jump function $j$ 
are the one associated with fractional Laplacian: $\mu(dz) = dz / {|z|^{N+\alpha}}$ 
(with $0 < \alpha < 2$) and $j(x,z)=z$. This class of operators appears in the context of stochastic control of jump processes and it is also worth pointing out that, at least to the best of our knowledge, general comparison results (for viscosity solutions) are only available for operators in the form \eqref{op:li}.  
\medskip

Many papers deal with H\"older estimates for fully non-linear elliptic equations. There are two kinds of approaches: for uniformly elliptic equations, one can use the powerful approach by Harnack inequalities which leads also to further regularity results; we refer the reader to Cabr\'e and Caffarelli \cite{CClivre} or Trudinger \cite{Tr1,Tr2} and references therein for results in this direction. A simpler method, more closely related to classical viscosity solutions theory, was introduced
by Ishii and Lions in \cite{il90}. Besides of its simplicity, it has the further advantage to provide results under weaker ellipticity assumptions and even for some degenerate equations; it was next used in \cite{gb91}, \cite{CY} where further regularity results are also proved and later in \cite{bs01,bd06}. As far as integro-differential elliptic equations are concerned, many papers were published in potential theory and equations are linear in most of these papers; moreover they rely on probabilistic techniques. See for instance \cite{bl02}. One of the first paper about H\"older estimates for integro-differential equations
with PDE techniques is probably \cite{silvestre06}. The author mainly deals with linear equations where singular measures
$\mu_x$ have a very general $x$-dependence, improving the previous literature. 
He is also able to deal with quite particular non-linear equations; more precisely, { he treats the case of equations like $F(\mathcal{I}[x,u],
\mathcal{J}[x,u])=0$ where $\mathcal{I}$ and $\mathcal{J}$ are two different non-local terms, and for ``strictly elliptic'' functions $F$ (in a suitable sense)}. 
\medskip

In the present paper, we deal with fully non-linear elliptic equations 
and we obtain local $C^{0,\alpha}$ regularity and estimates for a quite general class of 
integro-differential operators,
namely operators of the form \eqref{op:id} satisfying proper assumptions
(see \eqref{hyp:mu-1}-\eqref{hyp:mu-3}). 
Even if 
important operators \eqref{op:li} can be seen as special cases of the general one, 
we will give specific results with specific assumptions. {In other words, our second
main theorem is not a corollary of the first one}. 

Let us mention that, on the one hand, we cannot treat all 
the examples given in \cite{silvestre06} (in particular the examples of the Section 3.4
of \cite{silvestre06}) and, on the other hand, 
we can treat examples out of the scope of \cite{silvestre06};
indeed, we only assume that the measure $\mu_x$
is bounded at infinity (uniformly with respect to $x$), while Condition~(2.2) of \cite
{silvestre06} requires a (small) power of $|z|$ to be integrable at infinity. 
Moreover, we can handle much more general non-linear equations, in particular the important Bellman-Isaacs equations and we can also identify the critical H\"older exponent $\alpha$ of the solution. 
To be more precise, we are able to prove that the solution is $\alpha$-H\"older continuous 
for any $\alpha < \beta$ (and even $\alpha = \beta$ under stronger assumptions in the case $\beta <1$) 
where $\beta$ characterizes the singularity of the measure associated 
with the integral operator. 
\medskip

In order to treat a large class of non-linear elliptic equations, we decided
to present the main results by assuming that the non-linearity $F$ in \eqref{eq:pide}
satisfies a proper ellipticity-growth condition (see (H) in Section~\ref{sec:main}). 
Loosely speaking, 
this structure condition ensures that either the equation is strictly elliptic in 
the classical fully non-linear sense or it is strictly elliptic in a 
non-linear and non-local sense. Since this condition is rather technical, a whole section
is devoted to comments and examples (see Section~\ref{sec:examples}). 
\medskip

The techniques we use in the present paper do not seem to yield Lipschitz regularity
when $\beta \ge 1$ and we would like to investigate this question in a future work. 
Let us mention that we  proved \cite{bci07}
in quite a general framework that there exists a viscosity solution of the Dirichlet
problem without loss of boundary condition. The local H\"older estimates we obtain
in the present paper applies to the Dirichlet problem too 
and we will address naturally in a future work the question of 
boundary estimates.  We would like also to point out that the techniques we develop here can be readily applied to 
parabolic integro-differential equations. Finally, another possible interesting
application of these results is the study of the ergodicity of non-local equations.
\medskip

We would like to conclude this introduction by mentioning that after this work was
completed, we found out that Caffarelli and Silvestre \cite{csil} obtained regularity results for
a large class of nonlinear integrodifferential equations that are invariant under $x$-translations
and uniformly elliptic in a non-local way. In particular, they are able to get 
$C^{1,\alpha}$ regularity of solutions.
\medskip

The paper is organized as follows. In Section~\ref{sec:visc}, we recall the definition
of a viscosity solution of \eqref{eq:pide}. { In Section~\ref{sec:main}, we state two
main results:} the first one deals with general $x$-dependent L\'evy measure and the second
one deals with integro-differential operators under the L\'evy-It\^o form.
In Section~\ref{sec:examples}, we make comments on the main structure assumption on the 
non-linearity $F$ and we give examples of direct applications of our results. 
Section~\ref{sec:proofs} is devoted to proofs of the main results. 
\medskip

{
\noindent \textbf{Acknowledgements.}
The authors want to thank the referees for their very careful reading of the article and for valuable suggestions which lead to a real improvement of the presentation of our results. They are particularly grateful to one of the referees who provided a technical idea which allows them to slightly weaken the assumptions of their results.}

\medskip

\noindent \textbf{Notation.} The scalar product in $\R ^N$ is denoted by 
$x\cdot y$. A ball centered at $x$ of radius $r$ is denoted by $B(x,r)$. 
If $x=0$, we simply write $B_r$ and if $r=1$, we even write $B$. 
$\mathcal{S}^{N-1}$ denotes the unit sphere of $\R^N$. 

The transpose of a matrix $A$ is denoted $A^*$ and $||A||$ stands for the usual norm of $A$, 
namely $ ||A||:= \max_{|z| = 1}\, |Az| $. 
$\S$ is the space of $N \times N$, real symmetric matrices. 
We recall that $X \ge Y$ if all the eigenvalues of $X-Y$ are non-negative.   
If $X \in \S$ and $\eps \in (0,1)$ is such that all the eigenvalues of $X$ 
are strictly less than $1$ (resp. strictly greater than $-1$), we set 
$X^\eps = (I - \eps X)^{-1}X$ (resp. $X_\eps = (I + \eps X)^{-1}X$). 
These matrices are obtained  from $X$
by applying a sup-convolution procedure (resp. an inf-convolution procedure),
 namely, for any $\xi \in \R^N$
$$
X^\eps \xi \cdot \xi = \sup_{\zeta\in \R^N} \left\{X \zeta \cdot \zeta
 - \frac{|\xi-\zeta|^2}{\eps}  \right\} \quad , \quad
X_\eps \xi \cdot \xi = \inf_{\zeta\in \R^N} \left\{X \zeta \cdot \zeta
 + \frac{|\xi-\zeta|^2}{\eps}  \right\}.
$$

\section{Viscosity solutions for PIDE}
\label{sec:visc}

In this section, we recall the notion of degenerate ellipticity for 
non-linear non-local equations and the definition of viscosity solutions
for such equations. 

\subsection{Degenerate ellipticity}

Throughout the paper, the domain $\Omega$ is an open subset of
$\R^N$ and the non-linearity $F$ is a continuous function. We also assume that \eqref{eq:pide}
is degenerate elliptic. In this framework, this means that we make the following
\medskip

\noindent \textsc{Assumption (E).} For any $x \in \R^N$, $u \in \R$, 
$p \in \R^N$, $X, Y \in \S$, $l_1, l_2 \in \R$
$$ F(x,u,p,X,l_1) \leq F(x,u,p,Y,l_2)\quad\hbox{if } X\geq Y,\ l_1 \geq l_2\; .$$

\subsection{Non-Local operators}

In order to define viscosity solutions for \eqref{eq:pide}, we introduce
two associated operators $\mathcal{I}^{1,\delta}$ and
$\mathcal{I}^{2,\delta}$
\begin{eqnarray*}
\mathcal{I}^{1,\delta} [x,p,u] &=& \int_{|z|< \delta} [u
(x+z) - u  (x) - p \cdot z \un_B(z) ] \mu_x (dz),  \\
\mathcal{I}^{2,\delta} [x,p,u] &=& \int_{|z|\ge \delta} 
[u (x+z) - u  (x) - p \cdot z \un_B (z) ] \mu_x (dz). 
\end{eqnarray*}

In the case of L\'evy-It\^o operators~\eqref{op:li},  
$\mathcal{I}^{1,\delta}$ and $\mathcal{I}^{2,\delta}$ are defined 
as follows
\begin{eqnarray*}
\mathcal{I}^{1,\delta} [x,p,u] &=& \int_{|z|< \delta} [u
(x+j(x,z)) - u  (x) - p \cdot j(x,z) \un_B (z) ] \mu (dz),  \\
\mathcal{I}^{2,\delta} [x,p,u] &=& \int_{|z|\ge \delta} 
[u (x+j(x,z)) - u  (x) - p \cdot j(x,z) \un_B (z)] \mu (dz).
\end{eqnarray*}

\subsection{Definition}

We now recall the definition of a viscosity solutions for \eqref{eq:pide}. 
We assume that we are given a function $u$ defined on the whole space
$\R^N$. 
\begin{defi}[Viscosity solutions]\label{def:visc-sol}
An upper semi-continuous (usc in short) function $u: \R^N \to \R$ is a \emph{subsolution} of \eqref{eq:pide} if
for any test function $\phi \in C^2 (B(x,\delta))$ such
that $u -\phi$ attains a  maximum on $B(x,\delta)$ at 
$x \in \Omega$,
$$
F(x,u(x), D \phi (x), D^2 \phi (x), \mathcal{I}^{1,\delta} [x,p,\phi]
+ \mathcal{I}^{2,\delta} [x,p,u]) \ge 0
$$
where $p = D \phi (x)$. \\
A lower semi-continuous (lsc in short) function $u: \R^N \to \R$ is a \emph{supersolution} of \eqref{eq:pide} if
for any test function $\phi \in C^2 (\R^N)$ such
that $u -\phi$ attains a maximum on $B(x,r)$ at $x \in \Omega$ for some 
$r>0$,
$$
F(x,u(x), D \phi (x), D^2 \phi (x), \mathcal{I}^{1,r} [x,p,\phi]
+ \mathcal{I}^{2,r} [x,u]) \le 0,
$$
where $p = D \phi (x)$.\\
A continuous function $u:\R^N \to \R$ is a \emph{solution} of \eqref{eq:pide} 
if it is both a sub and a supersolution. 
\end{defi}
\begin{rem}
It is possible to construct solutions of \eqref{eq:pide} in the case where
$\Omega = \R^N$. If $\Omega \neq \R^N$, boundary conditions must be
imposed. For instance, as far as the Dirichlet problem is concerned, the
function $u$ can be prescribed outside of $\Omega$. See \cite{bci07} for
further details. 
\end{rem}
\begin{rem}
As remarked in \cite{bi07}, one can choose $r=0$ in the previous definition,
at least in the case of L\'evy-It\^o operators. See \cite{bi07} for further details. 
\end{rem}

\section{Main results}
\label{sec:main}

{ In this section, we state the two main results of this paper: the first one is concerned
with non-local operators of the form \eqref{op:id} (Theorem~\ref{thm:reg-int-mux}) and the second one 
provides other (distinct) results for L\'evy-It\^o operators~\eqref{op:li} (Theorem~\ref{thm:reg-int-li}). }

The two results rely on a structure condition imposed to the non-linearity $F$. 
In order to formulate it, we consider two functions $\Lambda_1, \Lambda_2: 
\overline \Omega \to [0, +\infty)$ such that $ \Lambda_1(x) + \Lambda_2 (x) 
\geq \Lambda_0 >0$ on $\overline \Omega$.
\medskip

\noindent
(H) {\em (Ellipticity-Growth condition)} For any $R >0$, there exist constants 
$k \geq 0$, $\tau, \theta \in (0,1]$,  a locally bounded function 
$\chi: \R^+ \times \R^+ \to \R^+$, a modulus of continuity $\omega_F: (0,+\infty) \to (0, +\infty)$, 
$\omega_F (0+) = 0$ and two constants $\eta, \bar\eps_0 >0$ such that for any $x,y \in \Omega$ with 
$|x-y|\leq \eta$, $u,v \in \R$ with $|u|, |v|  \leq R$, $p,q\in \R^N$, $|q| \leq R$, $l_1\leq l_2$, $\varpi \in (0,1/3)$, $L>0$, $\bar\eps\in (0, \bar\eps_0)$ and $\hat{a}\in \mathcal{S}^{N-1}$,
we have
\begin{eqnarray*}
F(y,u,p,Y,l_2) &-& F(x,v,p + q ,X ,l_1)  \\
&\leq&  \Lambda_1 (x) \bigg( \tr(X-Y) + \frac{\omega_F (|x-y|)}{\bar{\eps}}
 + |x-y|^\tau |p|^{2+\tau} + |p|^2 + \chi(L,\eta) \bigg)  \\
&&+ \displaystyle  \Lambda_2 (x) \bigg( (l_1-l_2) 
 +  \frac{|x-y|^{2\theta}}{\bar \eps} + |x-y|^\tau |p|^{k+\tau} + C_F |p|^{k} 
+ \chi(L,\eta)\bigg) \; ,
\end{eqnarray*}
if the matrices $X,Y$ satisfy
\begin{equation}\label{cond-xy}
-\frac{4}{\bar \eps} 
I \le
\left[\begin{array}{rr} X & 0\\ 0 & -Y \end{array}\right] 
 \le  \frac{2}{\bar \eps} \left[\begin{array}{rr}Z &- Z\\-Z & Z\end{array}\right] + 
L\left[\begin{array}{rr} I & 0\\ 0 & 0 \end{array}\right] \; , 
\end{equation}
where $I$ denotes the identity matrix and $Z = I - (1+\varpi) \hat{a} \otimes \hat{a}$.
\medskip

In the next subsection, we will make comments on this structure condition and 
give several examples. The general results are the following. 
\begin{thm}[H\"older continuity for general non-local operators]\label{thm:reg-int-mux}
Assume that the measures $\mu_x$ satisfy the following: there
exist $\beta \in (0,2)$, a constant ${\tilde C_\mu} >0$, a modulus of continuity $\omega_\mu: (0,+\infty) \to (0, +\infty)$, 
$\omega_\mu (0+) = 0$ and, for $\eta \in (0,1)$, a constant $C_\mu (\eta)>0$ such that
for any $x,y\in \Omega$, $d \in \mathcal{S}^{N-1}$, $\eta, \delta \in (0,1)$, 
\begin{equation}
\label{hyp:mu-1}
\int_B |z|^2 \mu_x (dz) + \int_{\R^N \setminus B} \mu_x (dz) \le {\tilde C_\mu }
\quad , \quad
\int_{\{ z : |z|\le \delta, |d \cdot z| \ge (1-\eta) |z| \}} |z|^2 \mu_x (dz) \ge C_\mu (\eta) \delta^{2-\beta} 
\end{equation}
\begin{equation}
\label{hyp:mu-2}
\int_{B_\delta} |z|^2 |\mu_x-\mu_y| (dz) \le
    \omega_\mu (|x-y|) \delta^{2-\beta} 
\end{equation}
{ 
\begin{equation}
\label{hyp:mu-3}
\int_{B \setminus B_\delta} |z| |\mu_{x} -\mu_{y}|(dz) 
\leq
\begin{cases}
\omega_\mu (|x-y|) \delta^{1-\beta} & \hbox{if  }\beta \neq 1\; ,
 \\
 \omega_\mu (|x-y|) |\ln \delta | & \hbox{if  }\beta = 1\; ,
\end{cases}
\end{equation}
}
with, if $\beta=1$,  $\omega_\mu (\cdot)$ such that
$\omega_\mu (r) | \ln r | \to 0$ as $r \to 0$.
Suppose also that the non-linearity $F$ satisfies (H) for some parameters
$k , \tau, \theta$.

{
(i)
 If
$$
\theta >  \frac12(2-\beta) \qquad \mbox{ and }\qquad
\left\{\begin{array}{ll} k = \beta & \mbox{ if } \beta >1, \\
k<\beta & \mbox{ if } \beta \leq 1,  \end{array}\right.
$$
then any bounded continuous viscosity solution $u$ of \eqref{eq:pide}
is locally $\alpha$-H\"older continuous for $\alpha$ small
enough. Precisely, $\alpha$ must satisfy: 
$\alpha < 1$ if $\beta \ge 1$ and $\alpha <\frac{\beta-k}{1-k}$ if
$\beta <1$.

\smallskip

(ii) If $\beta < 1$ and if we assume moreover that $C_F=0$ in (H) and $\tau > k (\beta^{-1}-1)$,
then $u$ is $\beta$-H\"older continuous.

\smallskip

Moreover, in both cases (i) and (ii), the $C^{0,\alpha}$ and $C^{0,\beta}$ estimates depend
on $||u||_\infty$,  $N$ (dimension), the constants $\tilde C_\mu$, $C_\mu (\eta)$ and the function $\omega_\mu$ 
appearing in \eqref{hyp:mu-1}-\eqref{hyp:mu-3},
on the constants and functions appearing in (H).
}
\end{thm}
We now turn to L\'evy-It\^o operators. 
\begin{thm}[H\"older continuity with L\'evy-It\^o operators]\label{thm:reg-int-li}
Assume that the function $j$ appearing in the definition of 
$\mathcal{I}_{LI}$ satisfies: there exist $c_0, C_0 >$ such that, for any $x \in \Omega$ and $z \in \R^N$, 
\begin{equation} \label{cond:j-0}
\left\{ \begin{array}{l}
c_0 |z| \le |j(x,z)| \le C_0 |z| \\
 |j(x,z) - j (y,z) | \le C_0 |z| |x-y|^{\tilde{\theta}}
\end{array}\right.
\end{equation}
with $\tilde{\theta} \in (0,1)$. Assume, in addition, that the measure $\mu$ satisfies: there exist $\beta \in (0,2)$, a constant ${\tilde C_\mu} >0$ and, for any $\eta \in (0,1)$, a constant $C_\mu (\eta)>0$ such that
for any $x\in \Omega$, $d \in \mathcal{S}^{N-1}$, $\eta, \delta \in (0,1)$, 
\begin{equation} \label{cond:mu1}
\int_B |j(x,z)|^2 \mu (dz) + \int_{\R^N \setminus B} \mu (dz) \le {\tilde C_\mu }
\quad , \quad
\int_{\mathcal{C}_{\delta,\eta}(x)} |j(x,z)|^2 \mu (dz) \ge C_\mu (\eta) \delta^{2-\beta} 
\end{equation} 
where $\mathcal{C}_{\delta,\eta} (x):=\{z:\; |j(x,z)|\le \delta, \; |d \cdot j(x,z)| \ge (1-\eta) |j(x,z)| \}$, { and that, moreover, for $\delta$ small enough
\begin{equation}
\label{cond:mu2}
\int_{B \setminus B_\delta} |z| \mu (dz) 
\leq
\begin{cases}
{\tilde C_\mu} \delta^{1-\beta} & \hbox{if  }\beta \neq 1\; ,
 \\
{\tilde C_\mu} |\ln \delta | & \hbox{if  }\beta = 1\; ,
\end{cases}
\end{equation}
}
Assume, finally, that the non-linearity $F$ satisfies (H) with parameters $k,\tau,\theta$.
If
$$
\theta,\tilde{\theta} >  \frac12(2-\beta) \qquad \mbox{ and }\qquad
\left\{\begin{array}{ll} k = \beta & \mbox{ if } \beta >1, \\
k<\beta & \mbox{ if } \beta \leq 1,  \end{array}\right.
$$
then any bounded continuous viscosity solution $u$ of \eqref{eq:pide} with
$\mathcal{I}[x,u]$ replaced with \eqref{op:li}
is locally $\alpha$-H\"older continuous for any 
$\alpha < \min (1,\beta)$. 

If, in addition, $C_F=0$ in (H) and $\tau > k (1-\beta)\beta^{-1}$, then $u$ is $\beta$-H\"older
continuous when $\beta <1$. 

Moreover, the  $C^{0,\alpha}$ estimate depends
on $||u||_\infty$, $N$ and on the constants and functions appearing in (H) and \eqref{cond:j-0}-\eqref{cond:mu2}.
\end{thm}
\begin{rem} It is worth pointing out that \eqref{cond:j-0}-\eqref{cond:mu1}-\eqref{cond:mu2} are the analogues  
of \eqref{hyp:mu-1}-\eqref{hyp:mu-2}-\eqref{hyp:mu-3} {but they do not imply them. It is easy to see
that the first line of \eqref{cond:j-0} together with \eqref{cond:mu1} imply \eqref{hyp:mu-1} but \eqref{hyp:mu-2}
and \eqref{hyp:mu-3} do not derive from the second line of \eqref{cond:j-0} and \eqref{cond:mu2}.}

Typically we have in mind the measures $\mu$ which satisy for some $C_\mu^\pm >0$ and $\beta \in (0,2)$ 
\begin{equation} \label{cond:mu}
\frac{C_\mu^- }{|z|^{N+\beta}} dz
\le \mu (dz) \le \frac{C_\mu^+}{|z|^{N+\beta}} dz. 
\end{equation} 
and functions $j(x,z)$ such that $z \mapsto j(x,z)$
has an inverse function $J(x,Z)$ and that there exist $c_0, C_0 >0 $ such that 
\begin{equation}\label{cond:j}
\left\{ \begin{array}{rl}
\forall (x,z) \in B_r(x_0,0), &
c_0 |z| \le |j(x,z)| \le C_0 |z|, \\
\forall (x,Z) \in B_R (x_0,0), &
c_0 |Z| \le |j(x,Z)| \le C_0 |Z| \\
\forall Z \in \R^N, & c_0 \le |\det D_z J (x_0,Z) | \\
\forall z \in \R^N, & |j(x,z) - j (y,z) | \le C_0 |z| |x-y|^{\tilde{\theta}}
\end{array}\right.
\end{equation}
and these are the properties we will use. We are in such a case if, for instance, for any $x$, $D_z j (x, z)$ exists for $|z|$ small enough, is continuous in $(x,z)$ and non-singular for $z=0$. Such a condition appears in \cite{bl02}. 
\end{rem}

\section{Comments and examples}
\label{sec:examples}

In this section, we make comments on assumptions of the main theorems 
and give  examples of applications. 
More precisely, we illustrate the different terms appearing in the structure
condition (H); we give examples of non-local operators of type~\eqref{op:id} and 
\eqref{op:li}; eventually, we give a regularity result that applies to the Bellman-Isaacs
equation. 

\subsection{Non-Linearities}

{ 
In this subsection, we illustrate the structure condition (H) we used in 
Theorems~\ref{thm:reg-int-mux} and \ref{thm:reg-int-li} and, to do so, we consider the model equation
\begin{equation}\label{eqn:example}
- \mathrm{tr} \, (A(x)D^2 u) - c(x)\mathcal{I} [x,u] + H(x,u,Du) = 0\quad\hbox{in  }\Omega\; ,
\end{equation}
where $A: \Omega\to \S $, $c:\Omega\to \R$ and $H:\Omega\times\R\times\R^N\to\R$ are continuous functions 
and $\mathcal{I} [x,u]$ is a non-local term of type \eqref{op:id} or \eqref{op:li}.

First, Equation \eqref{eqn:example} has to be degenerate elliptic and therefore we assume: for all 
$x\in \Omega$, $A(x)\geq 0$ and $c(x) \geq 0$. For $A$, we are even going to use the more restrictive 
assumption (but natural in the probabilistic framework) 

$$
\text{for all } x\in \Omega, \; A(x)=\sigma (x) \sigma^* (x) 
$$
where $\sigma$ is a continuous function which maps $\Omega$ into the space
of $N \times p$-matrices for some $p \le N$.

We come back to the structure condition (H). It combines two different terms: the first one
permits to handle equations that are strictly elliptic in the usual sense. 
The second one permits to handle non-local equations that are strictly elliptic
in a generalized (non-local) sense. Notice that imposing $\Lambda_1 (x) + \Lambda_2 (x) \ge \Lambda_0 >0$
means that, at each point $x \in \Omega$, the non-linearity is either strictly
elliptic in the classical (non-linear) sense or strictly elliptic in the generalized
(non-local) sense.

A typical situation is the following: we are given two open subsets 
$\mathcal{O}_1, \mathcal{O}_2$ of $\Omega$ such that 
$\mathcal{O}_1 \cup \mathcal{O}_2 = \Omega$ and the closure of $(\mathcal{O}_1)^c$
 is included in $\mathcal{O}_2$; moreover we know that $F$ satisfies (H) in 
$\mathcal{O}_1$ with $\Lambda_1 (x) \equiv 1, \Lambda_2 (x) \equiv 0$ while 
$F$ satisfies (H) in $\mathcal{O}_2$ with $\Lambda_1 (x) \equiv 0, \Lambda_2 (x) \equiv 1$.
 Then, if $\Lambda$ is a continuous function in $\Omega$ which equals $1$ on the 
closure of $(\mathcal{O}_1)^c$ and which has support included in $ \mathcal{O}_2$, 
then it is easy to check that $F$ satisfies (H) with 
$\Lambda_1 (x) \equiv 1-\Lambda (x), \Lambda_2 (x) \equiv \Lambda (x)$.

\smallskip

On Equation~\eqref{eqn:example}, the structure condition (H) means that we assume
$$ A(x) \geq \Lambda_1 (x) I \quad\hbox{and}\quad c(x) \geq \Lambda_2 (x) \quad\hbox{in  }\Omega\; .$$
Typically this means that the second-order $- \mathrm{tr} \, (A(x)D^2 u)$ is uniformly elliptic in $\mathcal{O}_1$, while there is no degeneracy in the non-local variable $l$ in $\mathcal{O}_2$. Of course, in conditions ``$\Lambda_1 (x) \equiv 1$ in $\mathcal{O}_1$'' or 
``$\Lambda_2 (x) \equiv 1$ in $\mathcal{O}_2$'', the ``$\equiv 1$'' may be replaced by ``$\equiv \Lambda$'' with $\Lambda >0$.

Besides of this ellipticity assumption, we have to assume that $A$ (or more precisely $\sigma$) satisfies suitable continuity assumptions in $\mathcal{O}_1$ and $\mathcal{O}_2$: this appears in (H) in the second subterms of the $\Lambda_1$, $\Lambda_2$-terms. To describe these assumptions, we recall a standard computation which appeared for the first time in Ishii~\cite{ishii89}: we assume that $\sigma$ is bounded uniformly continuous in $\Omega$ and we denote by $\omega_\sigma$ its modulus of continuity. We apply \eqref{cond-xy}     
to the vector $z = (z_1,z_2)$ with 
$ z_1 = \sigma (\xb) e$, $ z_2 = \sigma (\yb) e$ and an arbitrary 
$e \in \mathcal{S}^{N-1}$,
and get 
$$ 
\sigma^T (\xb)X \sigma (\xb)e\cdot e -  \sigma^T (\yb)Y \sigma
(\yb)e\cdot e \leq \frac{1}{\bar{\eps}} \omega_\sigma^2 (|\bar{x}-\bar{y}|)
+ L  ||\sigma||_\infty^2 
$$
(we used that $Z \le I$). 
Therefore
$$ 
\tr (A(\xb)X) - \tr(A(\yb)Y) \leq \frac{1}{\bar{\eps}} d \omega^2_\sigma(|\xb-\yb|)
+  L d \|\sigma\|_\infty^2. 
$$
Hence, choose $\omega_F (r)= d \omega^2_\sigma (r)$ and 
$\chi(L,\eta) = L d \|\sigma\|_\infty^2$. It is worth pointing out that, 
in the uniformly elliptic case (i.e. in $\mathcal{O}_1$), $A$ or $\sigma$ is only required to be continuous while, in $\mathcal{O}_2$, $\sigma$ has to be H\"older continuous (the same computations as above provide 
the $|x-y|^{2\theta}$--term if $\omega_\sigma (r) = Cr^\theta$) in order to take advantage of the ellipticity of the non-local term.

We now turn to the non-local term. Our main remark is the following: in our formulation, the non-local term $l$ is in fact $c(x)\mathcal{I} [x,u]$ and not only $\mathcal{I} [x,u]$. In that way, assumption (H) is obviously satisfied by Equation~\eqref{eqn:example}. On the other hand, in order to verify the assumptions of Theorem~\ref{thm:reg-int-mux}, one has to replace the measure $\mu_x (dz)$ with $\tilde \mu_x (dz):=c(x)\mu_x (dz)$ and check if conditions are satisfied by this new measure. In other words, in the case of \eqref{op:id}, we rewrite $c(x)\mathcal{I} [x,u]$ as
$$ \int_{\R^N} ( u (x+z ) - u (x)  - D u (x) \cdot z \un_B (z) ) \tilde \mu_x (dz) \; .$$
Therefore, the continuity assumptions on $\mu_x$  in Theorem~\ref{thm:reg-int-mux}
are indeed continuity assumptions on $\tilde \mu_x$ and, for this reason, they both contain continuity assumptions on $\mu_x$ and on $c$.

A priori a similar approach could be used for \eqref{op:li} but this really means that we use Theorem~\ref{thm:reg-int-mux} instead of Theorem~\ref{thm:reg-int-li} in this case.

It remains to consider the first order term $ H(x,u,Du)$ in Equation~\eqref{eqn:example}: this appears in (H) in  the third, fourth and (part of the) 
fifth subterms of $\Lambda_1$, $\Lambda_2$-terms. As it is classical for (local) equations, we have to require growth conditions with respect to $Du$: 
quadratic growth ($|p|^2$) when the equation is uniformly elliptic and a $|p|^k$-growth with $k$ depending on the measure when the strong ellipticity
 comes from the non-local term.

For example, in the case of the fractional Laplacian 
$(-\Delta)^\frac\beta2$, the natural growth turns out to be $k =  \beta$, even if 
Theorem~\ref{thm:reg-int-mux} shows that the case $\beta \leq 1$ is a little 
bit more particular. We refer the reader to Subsection~\ref{sec:non-loc} 
where an example of 
equation involving a fractional Laplacian is treated in details.

These growth conditions on the gradient have to be combined with the regularity of coefficients: we are able to treat gradient terms of the form $c(x) |D u (x)|^m$ with $m=2$ if $c$ is merely bounded and $m=2+\tau$ if $c$ is locally $\tau$-H\"older continuous. We leave details to the reader. }
\medskip

We say more about these assumptions in the next subsections. In particular,
we  treat equations that are not exactly of the form \eqref{eq:pide} but can
be handled with the same techniques (see Subsection~\ref{sec:bi}). 

\subsection{Singular measures}

The model singular measure is the L\'evy measure 
associated with the fractional Laplacian $(-\Delta)^{\frac\beta2}$. { In this case,
$d\mu_x(z)=d\mu(z)= dz/|z|^{N+\beta}$ with $0<\beta<2$. }

A second simple example of measure $\mu_x$ is $c(x,z) \mu (dz)$ with 
a L\'evy measure $\mu$ satisfying \eqref{cond:mu1}-\eqref{cond:mu2} with $j(x,z)=z$ and 
 $c(x,z)$ satisfying for any $x,y \in \Omega$ and $z \in B$
$$ |c(x,z) - c(y,z)|\le \omega (|x-y|) \quad \hbox{where  }\omega(t) \to 0\hbox{ when }t\downarrow 0,
$$ and, for any $x\in \Omega$, $z \in \R^N$, $0< \underline c \leq c(x,z) \leq \overline c$, for some constants $\underline c , \overline c$.
One can thus easily check \eqref{hyp:mu-1}, \eqref{hyp:mu-2}, \eqref{hyp:mu-3} 
and Theorem~\ref{thm:reg-int-mux} applies for suitable non-linearities $F$. 
L\'evy measures associated with tempered stable L\'evy processes
satisfy \eqref{cond:mu1}-\eqref{cond:mu2}. Indeed, in this case
$$
\mu (dz) = 
 \un_{(0,+\infty)} (z) \left( G^+ e^{- \lambda^+ |z|} \frac{dz}{|z|^{N+\alpha}} \right) 
+\un_{(-\infty,0)} (z) \left( G^- e^{- \lambda^- |z|} \frac{dz}{|z|^{N+\alpha}} \right). 
$$
These measures appear in financial modeling, see for instance \cite{cgmy}. 

\subsection{A non-local equation involving the fractional Laplacian\label{sec:non-loc}}

In order to illustrate further our results,
we next consider the following model equation{ 
$$
(-\Delta)^{\frac\beta2} u+b(x)|D u|^{k+\tau} + |D u|^r=0
$$
where $b\in C^{0,\tau}$, $0<\tau<1,\ 0<k,r<2$.  In this case, Condition $(H)$ is satisfied with
 $\Lambda_1=0$, $\Lambda_2(x)=1>0$,  $\theta>0$ is arbitrary and $\tau,k,r$ 
appear in the equation. It is easy to check that \eqref{hyp:mu-1} 
is satisfied with exponent $\beta$: first, it is a L\'evy measure and if 
$\mathcal{C}_{\delta,\eta}$ denotes $\{|z|< |\delta|,(d\cdot z)\ge (1-\eta)|z|\}$,
 then by homogeneity and symmetry of $\mu$,
 \begin{eqnarray*}
\int_{\mathcal{C}_{\delta,\eta}} |z|^2\mu(dz) 
&=&\frac{|\mathcal{C}_{\delta,\eta}|}{|B_\delta|}\int_{B_\delta}|z|^2 \mu(dz)
= \frac{|\mathcal{C}_{1,\eta}|}{|B_1|} \int_{B_\delta}|z|^2 \mu(dz)
= \frac{|\mathcal{C}_{1,\eta}|}{|B_1|} \int_{B_\delta}|z|^{2-\beta-N}dz\\
&=&C_\mu (\eta) \delta^{2-\beta}\,.
\end{eqnarray*}
Moreover, the other hypotheses we make on $\mu$ (namely \eqref{hyp:mu-2} and \eqref{hyp:mu-3})
are automatically 
satisfied since $\mu$ is independent of $x$ (in other words, one can choose 
$ \omega_\mu=0$).
Since $\beta<2$, we cannot allow here a quadratic growth for the gradient term; 
indeed Theorem 1 and 2 work only for $k,r\le\beta$ in the absence of local ellipticity. It is also worth pointing out that, if $r \leq 1$, we have $|p|^r -|p+q|^r \leq |q|^r \leq R^r$ and therefore, even if $\beta < 1$, any such $r$ works since the $R^r$ can be absorbed in the $\chi(L,\eta)$-term.}

\subsection{The Bellman-Isaacs equation}
\label{sec:bi}

Let us illustrate Theorem~\ref{thm:reg-int-li} by considering an important 
second-order non-linear elliptic integro-differential equations
appearing in the study of stochastic control of processes with 
jumps, namely the Bellman-Isaacs equation. Let us mention the
work of Jakobsen and Karlsen \cite{jk05} in the evolution case 
where the authors use completely different techniques.  
\begin{cor}\label{cor:isaacs}
Consider the following Bellman-Isaacs equation in $\R^N$
\begin{equation}\label{eq:isaacs}
c u + \sup_{\lambda \in \Lambda} \inf_{\gamma \in \Gamma} 
\bigg\{ - \frac12 \mathrm{Tr} (\sigma_\ab (x) \sigma^*_\ab (x) D^2 u ) - b_\ab
(x) \cdot D u - \mathcal{I}_{LI}^\ab [x,u]  - f_\ab (x)\bigg\} =0 
\end{equation} 
with $c \ge 0$ and where $\mathcal{I}_{LI}^\ab [x,u]$ is a family of L\'evy-It\^o operators
associated with a common L\'evy measure $\mu$ and  a family of functions $j_\ab (x,z)$. 
Assume that 
\begin{enumerate}[label=(\roman{*})]
\item $\mu$ verifies \eqref{cond:mu1}-\eqref{cond:mu2} with constants independent of $\lambda, \mu$,
\item there exist $c_0, C_0 >0$ and $\tilde{\theta} \in (0,1)$ such that for any $(\lambda, \gamma) \in \Lambda \times
\Gamma$, 
$j_\ab$ satisfies \eqref{cond:j-0},
\item $\sigma_\ab$, $b_\ab$ and $f_\ab$ satisfy for 
some $\theta \in (0,1)$ (and some constant $C_F >0$)
\begin{eqnarray*}
\forall \alpha, \beta, \quad \|\sigma_\ab \|_{0,\theta} 
+ \|b_\ab \|_{0,\theta}  +\|f_\ab \|_{0,\theta} \le C_F.
\end{eqnarray*}
\end{enumerate}
If $\theta,\tilde{\theta} > \frac12 (2-\beta)$, 
then any bounded viscosity solution $u$ of \eqref{eq:isaacs} is
$\alpha$-H\"older continuous for any 
$\alpha < 1$ if $\beta \ge 1$ 
and for $\alpha < \frac{\beta-k}{1-k}$ if $\beta <1$.
\end{cor}
\begin{proof}{ 
Remark that \eqref{eq:isaacs} is not exactly of the form
\eqref{eq:pide}. Nevertheless, the proof of Theorem~\ref{thm:reg-int-li}
we present below can be adapted to this framework. It is enough
to check that structure condition (H) is satisfied by the linear
equations 
$$
F_\ab (x,u,p,A,l)= c u  - \frac12 \mathrm{Tr} (\sigma_\ab (x) \sigma_\ab^* (x) A) - b_\ab
(x) \cdot p - l_\ab  - f_\ab (x)\; ,
$$
with constants and functions appearing in (H) independent on $\lambda$, $\gamma$.

Indeed, if it is the case, then we have just to use the standard inequality
$$ \sup_{\lambda}\inf_\gamma (\cdots)-\sup_{\lambda}\inf_\gamma (\cdots)
 \le \sup_{\lambda, \gamma} (\cdots - \cdots).
$$
and we can conclude in this more general case too. }
\end{proof}

\section{Proofs of Theorems~\ref{thm:reg-int-mux} and \ref{thm:reg-int-li}}
\label{sec:proofs}

We prove successively Theorems~\ref{thm:reg-int-mux} and \ref{thm:reg-int-li}. 
On one hand, {the proofs are very similar and we will skip details in the proof of the second
theorem when adapting arguments used in the proof of the first one}.
On the other hand, we need to use very precisely every parameter. This is the reason why
constants are computed from line to line and explicit formulae are given for each of them in order to use them later. 

\subsection{Proof of Theorem~\ref{thm:reg-int-mux}}

\begin{proof}[Proof of Theorem~\ref{thm:reg-int-mux}.~(i)]
Without loss of generality, we assume here that $C_F=1$ in the ellipticity-growth condition~(H). 
In order to prove the local 
H\"older continuity of $u$, we are going to show that, for any $x_0 \in \Omega$, 
there exists $L_2=L_2(x_0)$ such that, for some well chosen $\alpha \in (0,1)$ and
for $L_1=L_1(x_0)>0$ large enough, we have
$$
M = \sup_{x,y \in \R^N} \{ u (x) - u(y) -  \phi(x-y) - \Gamma (x)\} \le 0
$$
where $\phi (z) = {L_1} |z|^\alpha$ and  $\Gamma (x) = L_2 |x-x_0|^2$. 
We point out that the role of {the term $\Gamma (x)$} is to localize around $x_0$, 
while {the term $\phi (x-y)$} is concerned with the H\"older continuity. 
Proving such a result with a suitable control on $\alpha, L_1, L_2$ 
clearly implies the desired property.

If $\Omega \neq \R^N$, 
we first choose $L_2$ in order that $u (x) - u(y) -  \phi(x-y) - \Gamma (x) \leq 0$ 
if $x\notin \Omega$: to do so, we first choose
$$ L_2 \ge  \frac{8||u||_\infty}{[d(x_0,\partial\Omega)]^2}. $$
If $\Omega = \R^N$, $L_2$ is arbitrary. 
Then we argue by contradiction: we assume that $M>0$ and we are going 
to get a contradiction for $L_1$ large enough and for a suitable choice of $\alpha$.

If the supremum defining $M$ is attained at $(\bar{x},\bar{y})$, we then deduce from $M>0$ 
that $\bar{x} \neq \bar{y}$ and
\begin{equation}\label{penalisation}
|\bar{x} - \bar{y}| \le \left(\frac{2 \| u\|_\infty}{L_1}\right)^{\frac1\alpha}=: A, 
\qquad  
|\bar{x} -x_0| <  \sqrt{\frac{2\|u\|_\infty}{L_2}} =: R_2 \le \frac{d(x_0,\partial\Omega)}{2},
\qquad  
u(\bar{x}) > u (\bar{y}).
\end{equation}
{ If $L_1$ is large enough so that $\displaystyle  A< \frac{d(x_0,\partial\Omega)}{2}$, 
then we have $\bar{x} , \bar{y} \in \Omega$. }

Next, we pick some $\nu_0 \in (0,1)$ and we define 
$$
a= \bar{x}-\bar{y}, \qquad \qquad \eps = |a| \qquad \qquad 
\hat{a} = \frac{a}{|a|}, \qquad  \qquad \delta = \nu_0 \eps < \eps.
$$ 
{ First, $\nu_0$ will be chosen small enough but fixed (independent of $L_1$ and $\eps$). }

\medskip

From the study of the maximum point property for $(\bar{x},\bar{y})$, we also get
\begin{equation}\label{estim:omegau}
L_1 \eps^\alpha \le u (\xb) - u(\yb) \le \omega_u (\eps)
\end{equation}
where $\omega_u$ denotes the modulus of continuity of $u$ on $B(x_0,d(x_0,\partial\Omega)/2)$. 
We will use this piece of information   below (see Step~4). 
Notice also that if $\chi(x,y)$ denotes $\phi(x-y)+\Gamma(x)$, then
$x \mapsto u(x) - \chi(x,\bar{y})$ (resp. $y \mapsto u(y) +
\chi(\bar{x},y)$) attains a global maximum (resp. minimum) at $\bar{x}$
(resp. $\bar{y}$)  with $\chi(\cdot,\bar{y})$ (resp. $-\chi(\bar{x},y)$)
of class $C^2$  on $B(\bar{x},\delta)$. In particular, we can use $\chi(\bar{x},\cdot)$ and
$-\chi(\cdot,\bar{y})$ as test-functions in
Definition~\ref{def:visc-sol} with any $\delta' < \delta$. 
\smallskip

The remaining of the proof is divided in four steps. We 
write down the viscosity inequalities and combine them (Step~1), 
we get suitable matrices inequalities from non-local Jensen-Ishii's lemma (Step~2), 
we estimate from above the difference of the non-local terms (Step~3) 
and we conclude (Step~4). 
\medskip

\noindent \textsc{Step 1: writing down viscosity inequalities.}
Let $p$ denote $D \phi (a)$ and $q$ denote $D \Gamma (\bar{x})$. 
We use Corollary~1 of \cite{bi07} which, for $\iota>0$ small enough, provides us 
with two matrices $X_\iota, Y_\iota \in \S$ 
such that, for any $\delta' \ll 1$ 
\begin{eqnarray*}
F(\bar{x},u(\bar{x}),p+q,X_\iota, \mathcal{I}^{1,\delta'} [\bar{x},p+q,\chi_\iota(\cdot,\bar{y})]
+\mathcal{I}^{2,\delta'}[\bar{x},p+q,u]  + o_\iota(1)) \le 0, \\
F(\bar{y},u(\bar{y}),p,Y_\iota, 
\mathcal{I}^{1,\delta'} [\bar{y},p,-\chi_\iota(\bar{x},\cdot)]+\mathcal{I}^{2,\delta'} [\bar{y},p,u]
+ o_\iota(1)) \ge 0
\end{eqnarray*}
(where $\chi_\iota$ is an approximation of $\chi$ by a localized inf-convolution, see
 \cite{bi07})
with  the following matrix inequality
\begin{equation} 
- \frac{1}{\iota} I  \le \left[\begin{array}{rr}X_\iota & 0 \\ 0 & -Y_\iota \end{array}\right] 
 \le \left[\begin{array}{rr}Z&-Z\\-Z&Z\end{array}\right] + 
2L_2 \left[\begin{array}{rr}I&0\\0&0\end{array}\right] + o_\iota(1),
\label{ineg-mat}
\end{equation}
where $Z= D^2 \phi (a) $. 

Our aim is to first let $\iota$ tend to $0$ in order to get rid of all the 
artificial $\iota$ dependences in these inequalities: in order to do so, but also in order 
to apply (H) which requires a two-side bound on the matrices, we are first going
 to build matrices $X,Y$ such that the above viscosity inequalities still hold 
if we replace $X_\iota, Y_\iota$ by $X,Y$ and such that the matrices $X,Y$ satisfy 
the required inequality in (H).

Then, if we set
$$
l_1:=\mathcal{I}^{1,\delta'} [\xb,p+q,\chi_\iota (\cdot , \yb)]+\mathcal{I}^{2,\delta '}
[\xb,p+q,u], \qquad l_2:= \mathcal{I}^{1,\delta'} 
[\yb,p,-\chi_\iota (\xb, \cdot )]+\mathcal{I}^{2,\delta'} [\yb,p,u]\; ,
$$
and if we subtract the viscosity inequalities, 
dropping all the $\iota$ dependences, we will have
\begin{equation}\label{vis}
 0 \leq F(\yb,u(\yb),p,Y, l_2) - F(\xb,u(\xb),p+q,X, l_1)\; .
\end{equation}
In order to get the desired contradiction, the rest of the proof consists 
in obtaining various estimates, and in particular on the differences $X-Y$ 
and $\mathcal{I}^{2,\delta'} [\xb,p+q,u] - \mathcal{I}^{2,\delta'} [\yb,p,u]$, 
in order to apply the ellipticity-growth condition (H) to show that the 
right-hand side of this inequality is strictly negative. We point out that, 
since we are going to let first $\delta ' $ tend to $0$, the terms 
$\mathcal{I}^{1,\delta'} [\xb,p+q,\chi_\iota (\cdot , \yb)]$, 
$\mathcal{I}^{1,\delta'} [\yb,p,-\chi_\iota (\xb, \cdot )]$ create no difficulty because 
they tend to $0$ with $\delta'$.

\medskip
\noindent \textsc{Step 2: building and estimating the matrices $X, Y$.} We follow here
 ideas introduced by Crandall and Ishii \cite{CI} to obtain these matrices, by using only 
the upper bounds on $X_\iota, Y_\iota$. 
We rewrite Inequality~\eqref{ineg-mat} as: for any $z_1,z_2 \in \R^N$, we have
$$ 
X_\iota z_1 \cdot  z_1 - Y_\iota z_2 \cdot  z_2 \leq Z
(z_1 -z_2)\cdot  (z_1 -z_2) + 2L_2 |z_1|^2 .
$$
We have dropped the $o_\iota (1)$-term in the right-hand side for the sake of simplicity 
since it plays no role.
In fact, we use the previous matrix inequality on the following form 
$$ 
(X_\iota-2L_2 I) z_1 \cdot  z_1 - Y_\iota z_2 \cdot  z_2 \leq Z
(z_1 -z_2)\cdot  (z_1 -z_2). 
$$

Next we have to compute $Z$, as well as, for the rest of the proof, the derivatives of 
$\phi$. It will be convenient to do the proof for $\phi (x)= \varphi (|x|)$ for a general
smooth function $\varphi: \R^+ \to \R$. We thus get for any $b \in \R^N$
\begin{eqnarray}\label{gen-comp}
D\phi (b) &=& \varphi'(|b|)\hat{b} \nonumber \\
D^2 \phi (b) &=& \frac{\varphi'(|b|)}{|b|} P_{b^\perp} + \varphi'' (|b|) \hat{b} \otimes \hat{b}
\end{eqnarray}
where $\hat{b} = b / |b|$ and $P_{b^\perp} = I - \hat{b} \otimes \hat{b}$ is the projection 
on the orthogonal space of $b$. 
Hence, if $\varphi (r) = L_1 r^\alpha$, we get
\begin{eqnarray}
\nonumber
D \phi (b) &=& L_1 \alpha |b|^{\alpha-2}b \\
\label{hess-prec}
D^2 \phi (b) &=& L_1  (\alpha |b|^{\alpha-2} I + \alpha (\alpha-2) 
|b|^{\alpha-4} b \otimes b)
= L_1 \alpha |b|^{\alpha-4} (|b|^2 I - (2-\alpha) b \otimes b)\\
D^2 \phi (b) & \le & L_1 \alpha |b|^{\alpha-2} I.
\label{hess-ineg}
\end{eqnarray}
We have
$$
Z = \frac1{\bar \eps} \bigg(I - (2-\alpha) \hat{a} 
\otimes \hat{a} \bigg) \quad \mbox{ with } \quad  \bar \eps=(L_1 \alpha \eps^{\alpha-2})^{-1} \; .
$$

Now we come back to the $X_\iota, Y_\iota$ inequality: we apply to this inequality 
a sup-convolution in both variables $z_1$ and $z_2$ with a parameter which is 
$\frac14 \bar \eps$. Noticing that this corresponds to an inf-convolution on the
 $(Y_\iota z_2 \cdot  z_2)$--term, we easily get,
 with the notation introduced at the end of the Introduction, 
$$ (X_\iota-2L_2 I)^{\frac14 \bar \eps} z_1 \cdot  z_1 
- (Y_\iota)_{\frac14 \bar \eps} z_2 \cdot  z_2 \leq Z^{\frac12 \bar \eps}
(z_1 -z_2)\cdot  (z_1 -z_2). $$
On the other hand,  tedious but easy, explicit computations { (which are provided in the Appendix)} yield
$Z^{\frac12 \bar \eps} = 2{\bar \eps}^{-1} (I - (1+\varpi) \hat{a}\otimes \hat{a})$
 with 
$$\varpi:= \frac{1-\alpha}{3-\alpha} >0.$$
Notice that $0 < \varpi < 1/3$.

If we set $X = (X_\iota-2L_2 I)^{\frac14 \bar \eps} + 2L_2 I$, $Y=(Y_\iota)_{\frac14 \bar \eps}$,
 then $X,Y$ satisfy \eqref{cond-xy} with $L=2L_2$ and since $X_\iota \leq X$ and  $Y \leq Y_\iota$,
 the viscosity inequalities still hold for $X$ and $Y$ because $F$ is degenerate elliptic.

From this new form of inequality~\eqref{ineg-mat}, we can obtain several {types} of estimates 
on $X$ and $Y$: first, choosing $z_2 = -z_1 =\hat{a}$, we get
$$
X \hat{a} \cdot  \hat{a} - Y \hat{a} \cdot  \hat{a} 
\leq - \frac{8\varpi}{\bar{\eps}} + O(L_2) = - 8 L_1 \alpha \varpi \eps^{\alpha-2} + O(L_2) \; .
$$
Next choosing $z_2 = z_1 = z$ with $z$ being orthogonal to $ \hat{a}$, we have
$$ X z\cdot  z - Y z \cdot  z \leq O(L_2)\; .$$
In particular, this yields 
\begin{equation}\label{estim:diff-tr}
\tr (X-Y) \leq - 8 L_1 \alpha \varpi \eps^{\alpha-2}  + O(L_2)\; .
\end{equation}
\medskip

\noindent \textsc{Step 3: estimates of the non-local terms}.
The difference of the non-local terms, denoted by $T_{nl}$, can be rewritten as follows { (we recall that $B= B(0,1)$)}
\begin{eqnarray*}
T_{nl} & = & 
\int_{|z|\ge \delta'} [u(\bar{x}+z) -u (\bar{x}) - (p+q) \cdot z \un_B (z)] \mu_{\bar{x}} (dz) 
- \int_{|z|\ge \delta'} [u(\bar{y}+z) -u (\bar{y}) - p \cdot z \un_B (z)] \mu_{\bar{y}} (dz).
\end{eqnarray*}
In order to estimate it, we part the domain of integration $\{ |z| \ge \delta'\}$ 
into three pieces $\R^N \setminus B$ which leads to the $T_1$ term below, $
\mathcal{C} = \{ z \in B_\delta: |z \cdot \hat{a}| \ge (1-\eta)|z| \} \subset B
$ which leads to the $T_2$ term below and $B \setminus \mathcal{C}$ which leads to 
the $T_3$ term below.

In order not to add further technicalities, we assume from now on
that $\delta'=0$; the reader can check that if $\delta'>0$, the estimates
we present below remain valid up to $o_{\delta'} (1)$. 

Therefore we have to estimate from above $T_{nl} = T_1+T_2+T_3+o_{\delta'} (1)$ with
\begin{eqnarray*}
T_1 &=& \int_{|z|\ge 1} [u(\bar{x}+z) -u (\bar{x})] \mu_{\bar{x}} (dz) 
- \int_{|z|\ge 1} [u(\bar{y}+z) -u (\bar{y})] \mu_{\bar{y}} (dz)\\
T_2 &= & 
\int_{ \mathcal{C}} [u(\bar{x}+z) -u (\bar{x}) - (p+q) \cdot z ] \mu_{\bar{x}} (dz) 
- \int_{ \mathcal{C}} [u(\bar{y}+z) -u (\bar{y}) - p \cdot z] \mu_{\bar{y}} (dz) \\
T_3 & = & \int_{ B \setminus \mathcal{C}} [u(\bar{x}+z) -u (\bar{x}) - (p+q) \cdot z] \mu_{\bar{x}} (dz) 
- \int_{ B \setminus \mathcal{C}} [u(\bar{y}+z) -u (\bar{y}) - p \cdot z] \mu_{\bar{y}} (dz)
\end{eqnarray*}
{For the reader's convenience, we recall that $p= D \phi (a)$ and $q = D \Gamma (\bar{x})$.}
\medskip

\noindent \textbf{Estimate of $T_1$.} 
Since $u$ is bounded and so are the measures $\mu_x$ away from the origin, we conclude
that $T_1$ is bounded, uniformly with respect to all the parameters we
introduced. More precisely,
\begin{equation}\label{estim0}
T_1 \le \; C_1 
\end{equation}
where 
\begin{equation}\label{defconstante:1}
C_1 = 4 \|u\|_\infty \; \sup_{x \in B(x_0,d(x_0,\partial \Omega)/2)} \mu_x (\R^N \setminus B).
\end{equation}
\medskip

\noindent \textbf{Estimate of $T_2$.} 
We estimate $T_2$
from above by using the definition of $M$. Indeed, its definition 
provides the following key inequality
\begin{equation}
\label{eq:cle}
\begin{array}{lll}
u (\bar{x} +d) -u (\bar{x}) - (p+q) \cdot d  & \le & 
u (\bar{y} +d') -u (\bar{y}) - p \cdot d' \\
&& + \bigg\{ \phi (a+ d-d') - \phi (a) - 
D \phi (a) \cdot (d-d') \bigg\} \\
& & +  \bigg\{ \Gamma (\bar{x}+d)-\Gamma(\bar{x}) -D \Gamma
(\bar{x}) \cdot d \bigg\}. 
\end{array} 
\end{equation}
We then use \eqref{eq:cle} with $d = z$ and $d'=0$  
(resp. with $d=0$ and $d'= z$). We obtain
\begin{eqnarray*}
T_2 & \le& \int_{ \mathcal{C}} [ \phi (a+z)
-\phi (a) -D \phi (a) \cdot z ] \mu_{\bar{x}} (dz) \\
&&+ \int_{ \mathcal{C}}
[\phi (a-z)-\phi (a) + D \phi (a) \cdot z ]  \mu_{\bar{y}} (dz)
 + \int_{ \mathcal{C}}
[\Gamma (\bar{x}+z)-\Gamma(\bar{x}) -D \Gamma
(\bar{x}) \cdot z ] \mu_{\bar{x}} (dz)
\end{eqnarray*}
{We now use a second-order Taylor expansion} in each integral. First, according to the form of $\Gamma$, we have
$$ \int_{ \mathcal{C}}
[\Gamma (\bar{x}+z)-\Gamma(\bar{x}) -D \Gamma
(\bar{x}) \cdot z ] \mu_{\bar{x}} (dz) = L_2 \int_{ \mathcal{C}}
|z|^2 \mu_{\bar{x}} (dz)= O(L_2)\; . $$
Next, for the two other terms and for a general $\phi (x) = \varphi(|x|)$, we recall \eqref{gen-comp} and we finally get
\begin{eqnarray*}
T_2  &\le& \frac{1}{2}\int_{\mathcal{C}} \sup_{t \in (-1,1)}
\left( \frac{\varphi'(|a+tz|)}{|a+tz|} P_{(a+tz)^\perp} z \cdot z 
+ \varphi''(|a+tz|) (\widehat{a+tz} \cdot z)^2 \right)
(\mu_{\bar{x}}+\mu_{\bar{y}}) (dz) + O(L_2).
\end{eqnarray*}
and using next that $\varphi(r) = L_1 r^\alpha$, we obtain
\begin{eqnarray*}
T_2 &\le& 
\frac{L_1\alpha}{2}  \int_{\mathcal{C}} \sup_{t \in (-1,1)}
|a + t z|^{\alpha -4} (|a+tz|^2 |z|^2 - (2-\alpha)
((a+t z) \cdot z )^2) (\mu_{\bar{x}}+\mu_{\bar{y}}) (dz) + O(L_2).
\end{eqnarray*}
{ We use the notation $b$ for $a+t z$ and we estimate $|b|$ and $b\cdot z$ for $z \in \mathcal{C}$ as follows}
\begin{eqnarray*}
|b| & \le & \eps + t |z| \le \eps + \delta \le (1+ \nu_0) \eps  \\
|b\cdot z| & = & |a \cdot z + t |z|^2 | \ge (1-\eta) \eps |z| -  \delta |z| 
\ge (1-\eta-\nu_0) \eps|z|  \\
|b|^2 |z|^2 - (2-\alpha) (b \cdot z)^2 &\le &
(1 + \nu_0)^2 \eps^2 |z|^2 - (2-\alpha) (1-\eta-\nu_0)^2 \eps^2 |z|^2 \\
& \le & ( (1+ \nu_0)^2 - (2-\alpha) (1-\eta- \nu_0)^2) \eps^2   |z|^2 
\end{eqnarray*}
{ where we choose $\eta, \nu_0$ small enough so that $1-\eta- \nu_0 >0$ and}
\begin{equation} \label{condparam:2} 
(2-\alpha) (1-\eta- \nu_0)^2 -(1+ \nu_0)^2 > 0.
\end{equation}
Hence
$$
T_2  \le  - L_1 C_4 \eps^{\alpha-2} \int_{\mathcal{C}}
|z|^2 (\mu_{\bar{x}}+\mu_{\bar{y}}) (dz) + O(L_2)
$$
with 
\begin{equation}\label{defconstante:4}
C_4 =  \frac\alpha2 \bigg((2-\alpha) (1-\eta- \nu_0)^2
-(1+  \nu_0)^2  \bigg) (1+ \nu_0)^{\alpha-4}. 
\end{equation}
Using \eqref{hyp:mu-1}, we thus obtain
\begin{equation}\label{estim5-}
T_2  \le  - L_1 C_4 C_\mu (\eta) \eps^{\alpha-2} \delta^{2-\beta} + O(L_2),
\end{equation}
and finally, 
\begin{equation}\label{estim5}
T_2  \le  - L_1 C_5 \eps^{\alpha-\beta} + O(L_2),
\end{equation}
with
\begin{equation}\label{defconstante:6}
C_5 = C_4 C_\mu (\eta) \nu_0^{2-\beta}. 
\end{equation}
\medskip

\noindent \textbf{Estimate of $T_3$.}
In order to estimate $T_3$ from above, it is convenient to introduce
\begin{eqnarray*}
\mathcal{U}_1 (z) & = & u(\bar{x}+z) - u (\bar{x}) - (p+q) \cdot z \\
\mathcal{U}_2 (z) & = & u(\bar{y}+z) - u (\bar{y}) - p \cdot z 
\end{eqnarray*}
and write 
$$
T_3 =\int_{B \setminus \mathcal{C}} \mathcal{U}_1 (z) \mu_{\bar{x}} (dz) 
-  \int_{B \setminus \mathcal{C}} \mathcal{U}_2 (z) \mu_{\bar{y}} (dz).  
$$
We first remark that \eqref{eq:cle} with, successively 
{$(d,d')=(z,z)$, $(d,d')=(z,0)$ and $(d,d')=(0,z)=0$} yields
\begin{equation}\label{estim-1de+}
\left\{\begin{array}{ll}
\mathcal{U}_1 (z) - \mathcal{U}_2 (z) &\le \bigg( \Gamma(\bar{x}+z)-\Gamma(\bar{x}) - D \Gamma(\bar{x}) \cdot z \bigg) \\
\mathcal{U}_1 (z) &\le \bigg(\phi(a+z) - \phi (a) -D  \phi (a) \cdot z\bigg) 
+ \bigg( \Gamma(\bar{x}+z)-\Gamma(\bar{x}) - D \Gamma(\bar{x}) \cdot z \bigg) \\
\mathcal{U}_2 (z) &\ge - \bigg(\phi(a-z) - \phi (a) + D \phi (a) \cdot z\bigg).
\end{array}\right.
\end{equation}
We next consider the signed measure $\mu (dz) = \mu_{\bar{x}} (dz) - 
\mu_{\bar{y}} (dz)$. It can be represented by using its Hahn-Jordan decomposition with two non-negative measures 
{ $\mu^\pm$: we write $\mu = \mu^+ -\mu^-$ where $\mu^+ , \mu^-$ are respectively the positive and negative part of the measure $\mu$}. We would like next to introduce
a measure $\min( \mu_{\bar{x}},\mu_{\bar{y}})$. {To make it precise, we use the Hahn decomposition of $\R^N$ with respect 
to $\mu$: if $K$ denotes the support of $\mu^+$, }
we define $\overline{\mu} = \un_K \mu_{\bar{y}}
+ (1-\un_K) \mu_{\bar{x}}$. 
We now rewrite $T_3$ with these measures. We use
$$
\mu_{\bar{x}} = \mu^+ + \bar{\mu} \quad \text{and} \quad \mu_{\bar{y}} = \mu^- + \bar{\mu}
$$
together with \eqref{estim-1de+} to get
\begin{eqnarray*}
T_3 & = & \int_{B \setminus \mathcal{C}} (\mathcal{U}_1(z) -\mathcal{U}_2(z)) \overline{\mu} (dz) \\
&&+ \int_{B \setminus \mathcal{C}} \mathcal{U}_1 (z) \mu^+ (dz)- 
\int_{B \setminus \mathcal{C}} \mathcal{U}_2 (z) \mu^- (dz) \\
& \le & \int_{B \setminus \mathcal{C}} [ \Gamma(\bar{x}+z)-\Gamma(\bar{x}) - D\Gamma(\bar{x}) \cdot z ] 
\overline{\mu} (dz) \\
&&+ \int_{B \setminus \mathcal{C}} [\phi(a+z)- \phi(a) -D \phi (a) \cdot z] \mu^+ (dz)
+ \int_{B \setminus \mathcal{C}}  [ \Gamma(\bar{x}+z)-\Gamma(\bar{x}) - D\Gamma(\bar{x}) \cdot z ]\mu^+ (dz)\\
&&  +
\int_{B \setminus \mathcal{C}}  [\phi(a-z)- \phi(a) + D \phi (a) \cdot z] \mu^- (dz).
\end{eqnarray*}
In order to estimate the right hand side of the previous inequality from above, we first remark that
$$ \int_{B \setminus \mathcal{C}} [ \Gamma(\bar{x}+z)-\Gamma(\bar{x}) - D\Gamma(\bar{x}) \cdot z ] 
(\overline{\mu} + \mu^+) (dz) \le L_2 \int_B |z|^2 \mu_{\bar{x}}(dz) \le {\tilde C_\mu} L_2 .$$
Next, for the two other terms, we split the integration domain into $B \setminus B_\delta$ and $B_\delta \setminus
\mathcal{C}$. { On $B_\delta \setminus
\mathcal{C}$, we use once again a second-order Taylor expansion for $\phi$ while, on $B \setminus B_\delta$, we use the concavity of the function $t \mapsto L_1 t^\alpha$ on $(0,+\infty)$ to obtain
\begin{eqnarray}
\phi(a + z')- \phi(a) - D \phi (a) \cdot z' & \leq &
L_1(|a|+|z'|)^\alpha - L_1 |a|^\alpha - D \phi (a) \cdot  z' \nonumber \\
& \leq & \alpha L_1 |a|^{\alpha-1}|z'| + |D \phi (a) \cdot z'|  \nonumber\\
& \leq & 2 \alpha L_1 |a|^{\alpha-1}|z'| \; . \label{conc-ineq}
\end{eqnarray}

Using \eqref{conc-ineq} for $z'=z$ and $-z$ and \eqref{hyp:mu-1}, we are lead to
\begin{eqnarray*}
T_3 & \le &{\tilde C_\mu} L_2 + L_1 \int_{B_\delta \setminus \mathcal{C}} (\eps-\delta)^{\alpha -2}
|z|^2 |\mu| (dz) + 2 \alpha L_1\eps^{\alpha-1} \int_{B \setminus B_\delta} 
 |z| |\mu| (dz) \\
& \le &{\tilde C_\mu} L_2 
+ (1-\nu_0)^{\alpha-2} L_1 \eps^{\alpha -2} \int_{B_\delta \setminus \mathcal{C}} |z|^2 
|\mu_{\bar x}-\mu_{\bar y}| (dz) \\
& & +  2\alpha L_1 \eps^{\alpha-1}
\int_{B \setminus B_\delta} |z| |\mu_{\bar x} -\mu_{\bar{y}}|(dz). 
\end{eqnarray*}
where $|\mu| = \mu^+ + \mu^-=|\mu_{\bar x} -\mu_{\bar{y}}|$.

We now use \eqref{hyp:mu-2} and \eqref{hyp:mu-3}. It is convenient to introduce
$$
\psi_\beta (r) = \left\{\begin{array}{ll}
r^{1 -\beta} & \mbox{ if } \beta \neq 1 \\
| \ln r| & \mbox{if } \beta = 1. 
\end{array}\right. 
$$
We finally get
\begin{equation}
\label{estim-4-mux}
T_3 \le
{\tilde C_\mu} L_2 
+ \nu_0^{\alpha-\beta} L_1 \omega_\mu (\eps) \eps^{\alpha -\beta} 
+  2 \alpha L_1 \omega_\mu (\eps) \eps^{\alpha -1} \psi_\beta (\nu_0 \eps) 
\end{equation}
for $\nu_0$ small enough. 
We use here that $\nu_0^{\alpha-\beta}$ controls $(1-\nu_0)^{\alpha-2}\nu_0^{2-\beta}$ 
from above. 

\medskip

\noindent \textbf{Final estimate for $T_{nl}$.}
Combining \eqref{estim0},\eqref{estim5}, \eqref{estim-4-mux}, the final estimates are~: for $\beta \neq 1$
\begin{equation}\label{estim-final-tl-mux}
T_{nl} \le -L_1 C_5 \eps^{\alpha -\beta}  +  o (\eps^{\alpha-\beta}) + O (L_2), 
\end{equation}
and for $\beta =1$
$$
T_{nl} \le - L_1 C_5 \eps^{\alpha-1}
+  \nu_0^{\alpha-1} L_1 \omega_\mu (\eps) \eps^{\alpha-1} 
+  2\alpha L_1 \omega_\mu (\eps) \eps^{\alpha -1} |\ln (\nu_0 \eps)| + O (L_2).   
$$
We see that if $\omega_\mu (r)$ satisfies $\omega_\mu (r) |\ln r| \to 0$ as
$r \to 0$ (which is the case when $\beta = 1$), \eqref{estim-final-tl-mux} still holds true in this case. 

\medskip

\noindent \textsc{Step 4: conclusion.} For all $\alpha <  1$, we deduce
from \eqref{estim-final-tl-mux} that $T_{nl} \le 0$ if $L_1$ is large enough. 
Using inequality~\eqref{vis} together with (H) with $L=R=O(L_2)$ 
and Estimates~\eqref{estim:diff-tr} \&
\eqref{estim-final-tl-mux}, 
and recalling that $\bar \eps:= (L_1 \alpha \eps^{\alpha -2})^{-1}$,
we are thus lead to 
\begin{equation}\label{eq:lambda12} 
0 \le \Lambda_1 (\xb) A_1 + \Lambda_2 (\xb) A_2 
\end{equation}
with  
\begin{eqnarray*}
A_1 & =  &   L_1 \alpha \eps^{\alpha -2} \bigg[ 
- 8\varpi  + \omega_F (\eps) \bigg]
+ \eps^\tau (\alpha L_1 \eps^{\alpha-1})^{2+\tau} + (\alpha L_1 \eps^{\alpha-1})^2 
+ \tilde \chi(L_2)  \\
A_2 &= &  
 \bigg[ -L_1 C_5 \eps^{\alpha -\beta}  +  o (\eps^{\alpha-\beta}) \bigg]
 + \frac{\eps^{2\theta}}{\bar \eps} + \eps^\tau (\alpha L_1 \eps^{\alpha-1})^{k+\tau} 
+ (\alpha L_1 \eps^{\alpha-1})^{k} + \tilde \chi(L_2) 
\end{eqnarray*}
where we have gathered in the $\tilde \chi(L_2)$-term
the terms that either depend on $L_2$ or are bounded. 
We use the assumption $\Lambda_1 + \Lambda_2 \ge \Lambda_0 >0$ by rewriting \eqref{eq:lambda12}  as follows
$$
0 \le \Lambda_0 \max ( A_1,  A_2) \; .
$$
To get the desired contradiction and to obtain the $C^{0,\alpha}$-estimate, it is enough to prove that 
$A_i < 0$ for $i=1,2$ and $L_1$ large enough, and to control the size of such $L_1$.
 
\medskip

As far as $A_1$ is concerned and as soon as $\alpha < 1$, we can ensure that $- 8\varpi  + \omega_F (\eps) \leq - 4 \varpi $ if $L_1$ is large enough; using \eqref{penalisation}, it is clear that, in order to do it, the size of $L_1$ depends only on $||u||_\infty$, $\alpha$ and $d(x_0, \partial \Omega)$. This yields an estimate of the type
$$ 
A_1 \le 
  L_1 \eps^{\alpha -2}\bigg[
- 4\varpi   \alpha +  \alpha^{2+\tau} (L_1 \eps^{\alpha})^{1+\tau} + \alpha^2 L_1 \eps^{\alpha} \bigg]
 + \tilde \chi(L_2) \; .
$$
Now there are two ways to estimate $L_1 \eps^{\alpha}$: either to use the first part of inequality \eqref{estim:omegau} which yields the estimate $L_1 \eps^{\alpha} \leq 2||u||_\infty$, or to use the second part and the estimate of $L_1 \eps^{\alpha} $ through the modulus of continuity of $u$. In the sequel, the strategy of the proof consists in proving the result for $\alpha$ small enough by using the first estimate of $L_1 \eps^{\alpha}$ and then to use this first step (which provides a modulus of continuity of $u$) to prove it for all $\alpha$ by using the second estimate of $L_1 \eps^{\alpha}$.

Using $L_1 \eps^{\alpha} \leq 2||u||_\infty$ in the above inequality yields
$$ 
A_1 \le 
  L_1 \eps^{\alpha -2}\bigg[
- 4\varpi   \alpha +  \alpha^{2+\tau} (2||u||_\infty)^{1+\tau} + 2 \alpha^2 ||u||_\infty \bigg]
 + \tilde \chi(L_2) \; .
$$
For $\alpha$ small enough (depending only on $||u||_\infty$ and $\tau$), the bracket is less than $- 2\varpi   \alpha <0$  and, recalling \eqref{penalisation}, it is clear that the right-hand side is (strictly) negative if $L_1$ is large enough (depending on $\alpha$, $\varpi $,  $ \tilde \chi(L_2)$). Hence, we get the desired inequality:  $A_1 <0$.

\medskip

For the $A_2$-term, we first write
\begin{eqnarray*}
A_2 & =&  -L_1 C_5 \eps^{\alpha -\beta}  +  o (\eps^{\alpha-\beta})
 + \eps^{2\theta} L_1 \eps^{\alpha-2} + \eps^\tau (\alpha L_1 \eps^{\alpha-1})^{k+\tau}
+ (L_1 \alpha  \eps^{\alpha-1})^{k} + \tilde \chi(L_2) 
\\
& =&  L_1 \eps^{\alpha-\beta} 
\bigg[ -C_5   +  o_\eps (1) + \eps^{2\theta-2 +\beta} 
+ \alpha^{k+\tau}\eps^{\beta-k} (L_1 \eps^{\alpha})^{k+\tau -1}
+ \alpha^{k} \eps^{\beta-k}(L_1 \eps^\alpha)^{k-1} \bigg] 
+ \tilde \chi(L_2) \\
& =&  L_1 \eps^{\alpha-\beta} 
\bigg[ -C_5   +  o_\eps (1) +  \alpha^{k+\tau}\eps^{\beta-k} (L_1 \eps^{\alpha})^{k+\tau -1}
+ \alpha^{k} \eps^{\beta-k}(L_1 \eps^\alpha)^{k-1} \bigg] 
+ \tilde \chi(L_2)
\end{eqnarray*}
where $C_5$ is given by \eqref{defconstante:6}; we also used $2 \theta +\beta-2 >0$. }

The key difference with $A_1$ is the fact that the exponents of the term $L_1 \eps^{\alpha}$ can be 
non-positive and we have to argue differently if it is indeed the case.

We have the following cases.

\noindent $\bullet$ If $\beta > 1$, since $C_5 \geq \alpha C'_5$ for some constant $C'_5$ independent of $\alpha$ 
(at least for $\alpha \leq 1/2$), then one can argue as for $A_1$ with
$k=\beta$ since $k+\tau>1$ and $k=\beta > 1$ 
and obtain the $C^{0, \alpha}$ regularity and estimates for $\alpha$
small enough.

\noindent $\bullet$ If $\beta \leq 1$, then we cannot use this argument anymore since $k$ must satisfy (at least) 
$k\leq 1$.
In order to conclude, it is enough to ensure
\begin{equation}\label{cond:a2}
\eps^{\beta-k} (L_1 \eps^{\alpha})^{k+\tau -1} = o_\eps (1) \qquad 
\mbox{and} 
\qquad
\eps^{\beta-k}(L_1 \eps^\alpha)^{k-1} = o_\eps (1). 
\end{equation}
Writing $\eps^{\beta-k}(L_1 \eps^\alpha)^{k-1} = L_1^{k-1}\eps^{\beta-k
  + \alpha (k-1)}$, we see that this term is $o_\eps (1)$ if $k < 1$  
and $\beta-k + \alpha (k-1) \geq 0$; notice that we do not know how to
compare, in full generality, $L_1$-terms and $\eps$-terms. The second
condition implies that $k < \beta$.  
In the same way, for the other term, either $k+\tau -1 \geq 0$ and the condition $\beta >k$ is sufficient or 
$k+\tau -1 < 0$ and we are lead to $\alpha \leq(\beta-k)/(1-k-\tau)$. Gathering all these informations yields 
the conditions
$$ 1 \geq \beta > k \quad \hbox{and} \quad  \alpha \leq \frac{\beta-k}{1-k}\; .$$

At this point, putting together the informations on $A_1$ and $A_2$, we have shown that $u$ is locally 
in $C^{0,\bar \alpha}$ for $\bar \alpha$ small enough (depending only on the data and the $L^\infty$-norm of $u$) 
and we have an estimate of the local $C^{0,\bar \alpha}$-norm of $u$. In order to conclude the proof, we need 
to come back to the estimate on $A_1$ and $A_2$ (in the case when $\beta>1$) and to estimate the terms 
$L_1 \eps^{\alpha}$ using \eqref{estim:omegau} and the local 
$C^{0,\bar \alpha}$-modulus of continuity. This easily yields the full result and the { proof of point (i) is complete}.

{
In order to prove the second part of the theorem, we now choose $\alpha = \beta$ and we only need to adapt the final
step (Step~4) in the previous proof. We proceed as before by writing \eqref{eq:lambda12}
with $A_1$ unchanged and, since $C_F=0$, $A_2$ is given by
$$
A_2 = L_1 
 \bigg[ - C_5
+ o(1)  \bigg]
 +  L_1 \eps^{\beta -2 + 2\theta} + \eps^\tau (\beta L_1 \eps^{\beta-1})^{k+\tau} 
+ \tilde \chi(L_2)\; ,
$$
which we can rewrite as
$$
A_2 = L_1 
 \bigg[ - C_5
+ o(1)  \bigg]
 +  L_1 \eps^{\beta -2 + 2\theta} + \beta L_1 \eps^{\tau + (\beta-1)(k+\tau)} 
+ \tilde \chi(L_2)\; .
$$
At this stage of the proof, $L_2$ is fixed and $L_1$ can be chosen large enough in order to control
the term $\tilde\chi(L_2)$. Next we notice that $\beta -2 + 2\theta >0$. In order to conclude, it is enough to have $\tau + (\beta-1)(k+\tau)>0$, i.e. $\tau > k(\beta^{-1}-1)$.

 The proof of Theorem~\ref{thm:reg-int-mux} is now complete. }
\end{proof}

\subsection{Proof of  Theorem~\ref{thm:reg-int-li}}

\begin{proof}[Proof of Theorem~\ref{thm:reg-int-li}] 
This proof follows along the lines of the previous one, the only difference is the way of getting 
Estimate~\eqref{estim-final-tl-mux} for the non-local term in the new framework and under the new assumptions. Let us explain this point. 

In order to shed light on the fact that  $j$ has to be H\"older continuous with respect to $x$,  
we let $\omega (r)$ denote $C_0 r^{\tilde{\theta}}$ and we will see at the end
of the present proof that $\omega$ has to be chosen as a power law.
Precisely, in this case, 
\begin{eqnarray*}
T_{nl} & = & 
\int_{|z|\ge \delta'} [u(\bar{x}+j(\bar{x},z)) 
-u (\bar{x}) - (p+q) \cdot j(\bar{x},z) \un_B (z)] \mu (dz) \\
& & - \int_{|z|\ge \delta'} [u(\bar{y}+j(\bar{y},z)) 
-u (\bar{y}) - p \cdot j(\bar{y},z) \un_B (z)] \mu (dz).
\end{eqnarray*}
We then write $T_{nl} = T_1+T_2+T_3$ with
\begin{eqnarray*}
T_1 &=& \int_{|z|\ge 1} [u(\bar{x}+j(\bar{x},z)) -u (\bar{x})] \mu (dz) 
- \int_{|z|\ge 1} [u(\bar{y}+j(\bar{y},z)) -u (\bar{y})] \mu (dz)\\
T_2 &= & 
\int_{z \in \mathcal{C}} [u(\bar{x}+j(\bar{x},z)) 
-u (\bar{x}) - (p+q) \cdot j(\bar{x},z)  ] \mu (dz) \\
&&- \int_{z \in \mathcal{C}} [u(\bar{y}+j(\bar{y},z)) 
-u (\bar{y}) - p \cdot j(\bar{y},z)] \mu (dz) \\
T_3 & = & \int_{z \in B, z \notin \mathcal{C}} [\dots] \mu_{\bar{x}} (dz) 
- \int_{z \in B, z \notin \mathcal{C}} [\dots] \mu_{\bar{y}} (dz)
\end{eqnarray*}

\noindent where $\mathcal{C}$ is defined in the following way
$$  \mathcal{C}:= \{z: \; |j(\frac{\xb + \yb}{2},z)| \leq \frac{\delta}{2} \; \hbox{and} \; |j(\frac{\xb + \yb}{2},z)\cdot \hat a | \geq (1-\frac{\eta}{2})|j(\frac{\xb + \yb}{2},z)|\} = \mathcal{C}_{\delta/2,\eta/2}(\frac{\xb + \yb}{2})\; ,$$
where the notation $\mathcal{C}_{\delta,\eta}(x)$  is defined in the statement of Theorem~\ref{thm:reg-int-li}.
Roughly speaking, $ \mathcal{C}$ is the analogue of the cone used in the proof of Theorem~\ref{thm:reg-int-mux} where we have replaced $z$ by $j(\frac{\xb + \yb}{2},z)$. Notice that, because of \eqref{cond:j-0}, $ \mathcal{C}\subset B$ if $\delta$ is small enough.

We have chosen such a set $ \mathcal{C}$ for the following reason: if $L_1$ is large enough (or equivalently $\varepsilon$ or $\delta$ are small enough)
\begin{equation}\label{incimp}
 \mathcal{C} \subset \mathcal{C}_{\delta,\eta}(\xb) \cap \mathcal{C}_{\delta,\eta}(\yb) \; ,
 \end{equation}
which means that both $j(\xb,z)$ and $j(\yb,z)$ are in the ``good'' cones.

To check these properties, we write
\begin{eqnarray*}
|j(\xb,z) \cdot \hat a| & \ge & |j(\frac{\xb + \yb}{2},z)\cdot \hat a| - |j(\xb,z)-j(\frac{\xb + \yb}{2},z)| 
\ge (1-\eta/2) |j(\frac{\xb + \yb}{2},z)| - |z| \omega (|\xb-\frac{\xb + \yb}{2}|) \\
&\ge&   (1-\eta/2) |j(\xb,z)| - (2-\eta/2) |z| \omega (\varepsilon/2) \ge
(1 -\eta/2 -(2-\eta)C_0^{-1} \omega (\varepsilon/2) | j(\bar{x},z)| \\
& \ge & (1-\eta) |j(\bar{x},z)|
\end{eqnarray*}
and 
$$
|j(\bar{x},z)| \le |j(\frac{\xb + \yb}{2},z)| + \omega (\varepsilon/2) |z| \le \frac\delta2 + \frac\delta{2c_0}\omega (\varepsilon/2)
\le \delta
$$
for $L_1$ such that 
$$
\omega (\varepsilon/2) \le \min \left( \frac{\eta c_0}{4 -\eta} ,c_0 \right) =
\min(C_0,c_0) = c_0
$$
for $\eta < 1$.

\medskip

\noindent \textbf{Estimate of $T_1$.} We remark that, thanks to the properties of $j$ and $\mu$, \eqref{estim0} still holds
true. 
\medskip

\noindent \textbf{Estimate of $T_2$.}
One can check that \eqref{estim5-} and \eqref{estim5} still hold true. 
Indeed, we use \eqref{cond:mu1}-\eqref{cond:mu2} in order to get \eqref{hyp:mu-1}.
More precisely, using \eqref{cond:mu1}-\eqref{cond:mu2} and recalling the computation 
we made in Subsection~\ref{sec:non-loc}, we obtain
$$
\forall \hat{a}, \qquad   \int_{\mathcal{C}}
|z|^2 \mu (dz) \ge  C_\mu \delta^{2-\beta}
$$
where $C_\mu = C(c_0,C_0,C_\mu^-,\eta,d,\beta)$. Notice that \eqref{condparam:2} is slightly modified and so is $C_4$ 
(and consequently $C_5$). 
\medskip

\noindent \textbf{Estimate of $T_3$.} 
In the case of L\'evy-It\^o operators, we estimate $T_3$ as follows.
By using \eqref{eq:cle} with $d=j(\bar{x},z)$ and $d'=j(\bar{y},z)$, we get
\begin{eqnarray*}
 T_3 = \int_{z \in B, z \notin \mathcal{C}} [ u(\bar{x}+j(\bar{x},z)) -u (\bar{x}) 
- (p+q) \cdot j(\bar{x},z) \\
- u(\bar{y}+j(\bar{y},z)) + u (\bar{y}) + p \cdot
j(\bar{y},z)  ] \mu (dz) \\
\le T_3^1 + T_3^2 
\end{eqnarray*}
with 
\begin{eqnarray*}
T_3^1 & =& \int_{z \in B, z \notin \mathcal{C}} 
[\Gamma(\bar{x}+j(\bar{x},z)) -\Gamma (\bar{x}) -
q \cdot j(\xb,z)] \mu (dz)\\
T_3^2 & =&  \int_{z \in B, z \notin \mathcal{C}}  
[\phi(a+\Delta(z)) -\phi(a) -p \cdot \Delta(z)] \mu (dz)
\end{eqnarray*}
where $\Delta(z)  =  j(\bar{x},z) -j (\bar{y},z)$. 
Let us first estimate $T_3^1$ as follows. 
$$
T_3^1 \le \frac12 \int_B
\quad \sup_{t \in (0,1)} \bigg[ D^2 \Gamma (\bar{x}+t j
(\bar{x},z)) j (\bar{x},z) \cdot j (\bar{x}, z) \bigg] \mu (dz) 
$$
and we deduce that
\begin{equation}
T_3^1 \le  C_2 L_2  \label{estim2}
\end{equation}
where
\begin{equation}\label{defconstante:2}
C_2 =  C_0^2 \int_{|z|\le 1} |z|^2 \mu (dz).
\end{equation}
\medskip

We now turn to $T_3^2$ and we write $T_3^2 = T_3^{2,1}+T_3^{2,2}$ with
\begin{eqnarray*}
T_3^{2,1} & =&  \int_{z \in B, z \notin \mathcal{C}, 
|\Delta (z) | \ge \delta} [\phi (a+ \Delta (z))
- \phi (a) - D \phi (a) \cdot \Delta (z) ] \mu (dz), \\
T_3^{2,2} & =&  \int_{z \in B,z \notin \mathcal{C}, 
|\Delta (z) | \le \delta} [\phi (a+ \Delta (z))
- \phi (a) - D \phi (a) \cdot \Delta (z) ] \mu (dz). 
\end{eqnarray*}
We can now estimate $T_3^{2,1}$. In order to clarify computations, we write
$\omega$ for $\omega (\eps)$ in the following lines.

{ We use \eqref{conc-ineq} ; remarking that 
$$
\delta \le |\Delta (z)| \le \omega |z|\; ,
$$
we deduce
\begin{eqnarray*}
T_3^{2,1} & \le & \int_{\delta \omega^{-1} \le |z|\le 1} \,  2 \alpha L_1 \eps^{\alpha-1} |\Delta (z)|  \mu (dz) 
 \le  \int_{\delta \omega^{-1} \le |z|\le 1} 2 \alpha L_1
\eps^{\alpha-1} \omega |z| \mu (dz) \\
&\le & 2\alpha  L_1 \eps^{\alpha-1} \omega \psi_\beta (\delta \omega^{-1})  
\end{eqnarray*}
where, by \eqref{cond:mu2}
\begin{equation}\label{def:psi}
\psi_\beta (r) = \int_{r \le |z| \le 1} |z| \mu (dz) 
\le 
\begin{cases}
{\tilde C_\mu} r^{1-\beta} & \hbox{if  }\beta \neq 1\; ,
 \\
{\tilde C_\mu} |\ln r|& \hbox{if  }\beta = 1\; ,
\end{cases}
\end{equation}
In order to estimate $T_3^{2,2}$, we use a Taylor expansion
together with \eqref{hess-ineg} and the fact that $|a+ t\Delta(z)| \ge \eps-\delta > 0$. 
It comes
\begin{eqnarray*}
T_3^{2,2} &\le& 
\frac{L_1 \alpha}2 \int_{|z| \le 1} \sup_{t \in (0,1)}
|a+t \Delta (z)|^{\alpha-2} | \Delta (z)|^2 \mu (dz) \\
& \le & \frac{L_1}2 \int_B |z|^2 \mu (dz) (\eps -\delta)^{\alpha-2} \omega^2\\
& \le & \frac{L_1}2 \omega^2 \eps^{\alpha-2} (1 -\nu_0)^{\alpha-2}
\int_B |z|^2 \mu (dz)  \; .
\end{eqnarray*}
Gathering the estimates on $T_3^{2,1}$ and $T_3^{2,2}$, we obtain
\begin{equation}
T_3^2 \le  
2 \alpha  L_1 \eps^{\alpha-1} \omega \psi_\beta (\delta \omega^{-1}) 
+  L_1 C_3 \omega^2  \eps^{\alpha-2} 
\label{estim4}
\end{equation}
with 
\begin{equation}\label{defconstante:3}
C_3 =  \frac12 \int_B |z|^2 \mu (dz)  (1 -\nu_0)^{\alpha-2}  
\end{equation}
\bigskip

\noindent \textbf{Final estimate of $T_{nl}$.} 
Gathering estimates \eqref{estim0},  \eqref{estim5}, 
\eqref{estim2} and \eqref{estim4}, we finally obtain
\begin{equation} \label{estim6} 
T_{nl} \le  2\alpha L_1 \omega \psi_\beta (\nu_0 \omega^{-1}\eps) \eps^{\alpha-1} 
+ L_1 C_3 \omega^2 \eps^{\alpha-2} - L_1 C_5 \eps^{\alpha-\beta} 
+ O(L_2).
\end{equation}
Looking at the $\omega^2 \eps^{\alpha-2}$-term, we see that we need
to assume that $j$ is H\"older continuous with respect to $x$ and we thus replace
$\omega(\eps)$ with $C \eps^{\tilde{\theta}}$. Denoting $C'_3=CC_3$, we get 
$$
T_{nl} \le 2\alpha L_1 \eps^{\alpha+\tilde{\theta} -1} \psi_\beta (\nu_0 \eps^{1-\tilde{\theta}}) 
+ L_1 C'_3 \eps^{\alpha-2 + 2 \tilde{\theta}} - L_1 C_5 \eps^{\alpha-\beta} + O(L_2). 
$$
We claim that as in the case of Theorem~\ref{thm:reg-int-mux}, 
\eqref{estim-final-tl-mux} holds true under the assumptions of Theorem~\ref{thm:reg-int-li}.
To see this, we write
\begin{eqnarray*}
2 \alpha L_1 \eps^{\alpha+\tilde{\theta} -1} \psi_\beta (\nu_0 \eps^{1-\tilde{\theta}}) 
+ L_1 \eps^{\alpha-2 + 2 \tilde{\theta}}
&=& 
L_1 \eps^{\alpha-\beta} \bigg[
\eps^{\beta+\tilde{\theta} -1} \psi_\beta (\nu_0 \eps^{1-\tilde{\theta}})
+  C'_3 \eps^{\beta-2 + 2 \tilde{\theta}}\bigg] \\
&=&
L_1 \eps^{\alpha-\beta} \bigg[  \eps^{\beta+\tilde{\theta} -1} \psi_\beta (\nu_0 \eps^{1-\tilde{\theta}})
+ o_\eps(1)\bigg]
\end{eqnarray*} 
since $\tilde \theta > \frac12(2-\beta)$ and $\alpha < \beta$.
We next distinguish cases. 
\medskip

\noindent $\bullet$
If $\beta \neq 1$, then by \eqref{def:psi}, we get
\begin{eqnarray*}
\eps^{\beta+\tilde{\theta} -1} \psi_1 (\nu_0 \eps^{1-\tilde{\theta}}) 
= O (\eps^{\tilde\theta \beta}). 
\end{eqnarray*}
and we conclude in this case. 

\noindent $\bullet$
If $\beta =1$, 
\begin{eqnarray*}
\eps^{\beta+\tilde{\theta} -1} \psi_1 (\nu_0 \eps^{1-\tilde{\theta}}) 
= O ( \eps^{\tilde\theta} |\ln (\eps^{1-\tilde\theta})|) =o_\eps (1)
\end{eqnarray*}
and we can conclude in this case too. 

The proof of Theorem~\ref{thm:reg-int-li} is now complete. }
\end{proof}

{ 
\section*{Appendix}

We provide in this Appendix the explicit computation of $Z^{\frac12 \bar \eps}$ used in the proof of Theorem~\ref{thm:reg-int-mux}.

We have $Z = \frac1{\bar \eps} \bigg(I - (2-\alpha) \hat{a} 
\otimes \hat{a} \bigg)$ and 
$$Z^{\frac12 \bar \eps}q \cdot q =  \sup_{p \in \R^N} \left\{Z p \cdot p
 - \frac{2}{\bar\eps}|p-q|^2  \right\}= \frac{2}{\bar\eps}\sup_{p \in \R^N} \left\{\frac{\bar\eps}{2} Z p \cdot p
 - |p-q|^2  \right\} \; .$$
One checks easily that the (last) sup is achieved at a point $p$ such that $ \bar\eps Z p = 2(p-q)$, i.e.
\begin{equation}\label{lastone}
p - (2-\alpha) (\hat{a} \cdot p)  \hat{a} = 2(p-q)\; .
\end{equation}
Taking a scalar product with $\hat{a}$, we deduce 
$$ \hat{a} \cdot p = \frac{2}{3-\alpha} \hat{a} \cdot q\; ,$$
and inserting in the previous equality yields
\begin{equation}\label{lastonebis}
p -q = q -  \frac{2(2-\alpha)}{3-\alpha} (\hat{a} \cdot q) \hat{a}\; .
\end{equation}
Coming back to the value of the supremum and using \eqref{lastone}, we have
$$Z^{\frac12 \bar \eps}q \cdot q =  \frac{2}{\bar\eps} \left\{
(p-q) \cdot p - |p-q|^2\right\}= \frac{2}{\bar\eps} (p-q) \cdot q\; ,$$
and now using \eqref{lastonebis}, we finally obtain
$$
Z^{\frac12 \bar \eps}q \cdot q =  \frac{2}{\bar\eps} \left\{
|q|^2 - \frac{2(2-\alpha)}{3-\alpha} (\hat{a} \cdot q)^2 \right\} \; .$$
Therefore $1+\varpi = \frac{2(2-\alpha)}{3-\alpha}$ and $\varpi= \frac{1-\alpha}{3-\alpha}.$
}

\end{document}